\documentclass{amsart}

\usepackage{amsmath}
\usepackage{amsfonts}
\usepackage{amssymb,enumerate}
\usepackage{amsthm}
\usepackage[all]{xy}
\usepackage{hyperref}

\DeclareMathOperator{\hight}{ht}

\newcommand{\Spec}{\operatorname{Spec}}
\newcommand{\QSpec}{\operatorname{QSpec}}

\newcommand{\td}{\operatorname{tr.deg.}}
\newcommand{\GD}{\operatorname{GD}}
\renewcommand{\dim}{\operatorname{dim}}

\newcommand{\QMax}{\operatorname{QMax}}
\newcommand{\h}{\operatorname{ht}}

\newtheorem{thm}{Theorem}[section]
\newtheorem{cor}[thm]{Corollary}
\newtheorem{lem}[thm]{Lemma}
\newtheorem{prop}[thm]{Proposition}
\newtheorem{defn}[thm]{Definition}

\begin{document}

\bibliographystyle{amsplain}

\date{}

\author{Parviz Sahandi}

\address{Department of Mathematics, University of Tabriz, Tabriz,
Iran, and School of Mathematics, Institute for Research in
Fundamental Sciences (IPM), P.O. Box: 19395-5746, Tehran Iran}
\email{sahandi@ipm.ir}



\thanks{P. Sahandi was supported in part by a grant from IPM (No.
88130034)}

\title[Semistar dimension of polynomial rings]{Semistar dimension of polynomial rings and Pr\"{u}fer-like domains}

\begin{abstract} Let $D$ be an integral domain and $\star$ a semistar operation stable and of finite type on
it. In this paper we define the semistar dimension (inequality)
formula and discover their relations with $\star$-universally
catenarian domains and $\star$-stably strong S-domains. As an
application we give new characterizations of
$\star$-quasi-Pr\"{u}fer domains and UM$t$ domains in terms of
dimension inequality formula (and the notions of universally
catenarian domain, stably strong S-domain, strong S-domain, and
Jaffard domains). We also extend Arnold's formula to the setting of
semistar operations.
\end{abstract}

\maketitle

\section{Introduction}

\noindent All rings considered in this paper are (commutative
integral) domains (with 1); throughout, $D$ denotes a domain with
quotient field $K$. In \cite{OM}, Okabe and Matsuda introduced the
concept of a semistar operation. Let $D$ be an integral domain and
$\star$ a semistar operation on $D$.

In \cite{S} we defined and studied the $\widetilde{\star}$-Jaffard
domains and proved that every $\widetilde{\star}$-Noetherian and
P$\star$MD of finite $\widetilde{\star}$-dimension is a
$\widetilde{\star}$-Jaffard domain. In \cite{Sah} we defined and
studied two subclasses of $\widetilde{\star}$-Jaffard domains,
namely the $\widetilde{\star}$-stably strong S-domains and
$\widetilde{\star}$-universally catenarian domains and showed how
these notions permit studies of $\widetilde{\star}$-quasi-Pr\"{u}fer
domains in the spirit of earlier works on quasi-Pr\"{u}fer domains.
The next natural step is to seek a semistar analogue of dimension
(inequality) formula \cite{FHP}. In Section 2 of this paper we
define the $\widetilde{\star}$-dimension (inequality) formula and
show that each $\widetilde{\star}$-universally catenarian domain
satisfies the $\widetilde{\star}$-dimension formula and each
$\widetilde{\star}$-stably strong S-domain satisfies the
$\widetilde{\star}$-dimension inequality formula. In Section 3 we
give new characterizations of $\star$-quasi-Pr\"{u}fer domains and
UM$t$ domains in terms of the classical notions of dimension
inequality formula, universally catenarian domain, stably strong
S-domain, strong S-domain, and Jaffard domains. In the last section
we extend Arnold's formula to the setting of semistar operations
(see Theorem \ref{A}).

To facilitate the reading of the introduction and of the paper, we
first review some basic facts on semistar operations. Denote by
$\overline{\mathcal{F}}(D)$ the set of all nonzero $D$-submodules of
$K$, and by $\mathcal{F}(D)$ the set of all nonzero \emph{fractional
ideals} of $D$; i.e., $E\in\mathcal{F}(D)$ if
$E\in\overline{\mathcal{F}}(D)$ and there exists a nonzero element
$r\in D$ with $rE\subseteq D$. Let $f(D)$ be the set of all nonzero
finitely generated fractional ideals of $D$. Obviously,
$f(D)\subseteq\mathcal{F}(D)\subseteq\overline{\mathcal{F}}(D)$. As
in \cite{OM}, a {\it semistar operation on} $D$ is a map
$\star:\overline{\mathcal{F}}(D)\rightarrow\overline{\mathcal{F}}(D)$,
$E\mapsto E^{\star}$, such that, for all $x\in K$, $x\neq 0$, and
for all $E, F\in\overline{\mathcal{F}}(D)$, the following three
properties hold: $\star_1$: $(xE)^{\star}=xE^{\star}$; $\star_2$:
$E\subseteq F$ implies that $E^{\star}\subseteq F^{\star}$;
$\star_3$: $E\subseteq E^{\star}$ and
$E^{\star\star}:=(E^{\star})^{\star}=E^{\star}$. Let $\star$ be a
semistar operation on the domain $D$. For every
$E\in\overline{\mathcal{F}}(D)$, put $E^{\star_f}:=\bigcup
F^{\star}$, where the union is taken over all finitely generated
$F\in f(D)$ with $F\subseteq E$. It is easy to see that $\star_f$ is
a semistar operation on $D$, and ${\star_f}$ is called \emph{the
semistar operation of finite type associated to} $\star$. Note that
$(\star_f)_f=\star_f$. A semistar operation $\star$ is said to be of
\emph{finite type} if $\star=\star_f$; in particular ${\star_f}$ is
of finite type. We say that a nonzero ideal $I$ of $D$ is a
\emph{quasi-$\star$-ideal} of $D$, if $I^{\star}\cap D=I$; a
\emph{quasi-$\star$-prime} (ideal of $D$), if $I$ is a prime
quasi-$\star$-ideal of $D$; and a \emph{quasi-$\star$-maximal}
(ideal of $D$), if $I$ is maximal in the set of all proper
quasi-$\star$-ideals of $D$. Each quasi-$\star$-maximal ideal is a
prime ideal. It was shown in \cite[Lemma 4.20]{FH} that if
$D^{\star} \neq K$, then each proper quasi-$\star_f$-ideal of $D$ is
contained in a quasi-$\star_f$-maximal ideal of $D$. We denote by
$\QMax^{\star}(D)$ (resp., $\QSpec^{\star}(D)$) the set of all
quasi-$\star$-maximal ideals (resp., quasi-$\star$-prime ideals) of
$D$.

If $\Delta$ is a set of prime ideals of a domain $D$, then there is
an associated semistar operation on $D$, denoted by
$\star_{\Delta}$, defined as follows:
$$E^{\star_{\Delta}}:=\cap\{ED_P|P\in\Delta\}\text{, for each }E\in\overline{\mathcal{F}}(D).$$
If $\Delta=\emptyset$, let $E^{\star_{\Delta}}:=K$ for each
$E\in\overline{\mathcal{F}}(D)$. When
$\Delta:=\QMax^{\star_f}(D)$, we set
$\widetilde{\star}:=\star_{\Delta}$. It has become standard to
say that a semistar operation $\star$ is {\it stable} if $(E\cap
F)^{\star}=E^{\star}\cap F^{\star}$ for all $E$, $F\in
\overline{\mathcal{F}}(D)$. All spectral semistar operations are
stable \cite[Lemma 4.1(3)]{FH}. In particular, for any semistar
operation $\star$, we have that $\widetilde{\star}$ is a stable
semistar operation of finite type \cite[Corollary 3.9]{FH}.

The most widely studied (semi)star operations on $D$ have been
the identity $d_D$, and $v_D$, $t_D:=(v_D)_f$, and
$w_D:=\widetilde{v_D}$ operations, where
$E^{v_D}:=(E^{-1})^{-1}$, with $E^{-1}:=(D:E):=\{x\in
K|xE\subseteq D\}$.

For each quasi-$\star$-prime $P$ of $D$, the \textit{$\star$-height}
of $P$ (for short, $\star$-$\hight(P)$) is defined to be the
supremum of the lengths of the chains of quasi-$\star$-prime ideals
of $D$, between prime ideal $(0)$ (included) and $P$. Obviously, if
$\star=d_D$ is the identity (semi)star operation on $D$, then
$\star$-$\hight(P)=\hight(P)$, for each prime ideal $P$ of $D$. If
the set of quasi-$\star$-prime of $D$ is not empty, the
$\star$-dimension of $D$ is defined as follows:
$$\star\text{-}\dim(D):=\sup\{\star\text{-}\hight(P)|P\text{ is a
quasi-}\star\text{-prime of }D\}.$$ If the set of
quasi-$\star$-primes of $D$ is empty, then pose
$\star\text{-}\dim(D):=0$. Thus, if $\star =d_D$, then
$\star\text{-}\dim(D)= \dim(D)$, the usual (Krull) dimension of $D$.
It is known (see \cite[Lemma 2.11]{EFP}) that
$$
\widetilde{\star}\text{-}\dim(D)=\sup\{\hight (P) \mid P\text{ is a
quasi-}\widetilde{\star}\text{-prime ideal of } D\}.
$$

Let $\star$ be a semistar operation on a domain $D$. Recall from
\cite[Section 3]{EFP} that $D$ is said to be a
\emph{$\star$-Noetherian domain}, if $D$ satisfies the ascending
chain condition on quasi-$\star$-ideals. Also recall from
\cite{FJS} that, $D$ is called a \emph{Pr\"{u}fer
$\star$-multiplication domain} (for short, a P$\star$MD) if each
finitely generated ideal of $D$ is \emph{$\star_f$-invertible};
i.e., if $(II^{-1})^{\star_f}=D^{\star}$ for all $I\in f(D)$.
When $\star=v$, we recover the classical notion of P$v$MD; when
$\star=d_D$, the identity (semi)star operation, we recover the
notion of Pr\"{u}fer domain. Finally recall from \cite{CF} that
$D$ is said to be a \emph{$\star$-quasi-Pr\"{u}fer domain}, in
case, if $Q$ is a prime ideal in $D[X]$, and $Q\subseteq P[X]$,
for some $P\in\QSpec^{\star}(D)$, then $Q=(Q\cap D)[X]$. This
notion is the semistar analogue of the classical notion of the
quasi-Pr\"{u}fer domains. By \cite[Corollary 2.4]{CF}, $D$ is a
$\star_f$-quasi-Pr\"{u}fer domain if and only if $D$ is a
$\widetilde{\star}$-quasi-Pr\"{u}fer domain.

\section{The $\star$-dimension (inequality) formula}

We begin with the following definition. Recall that if $D\subseteq
T$ are domains, then $\td_D(T)$ is defined as the
\emph{transcendence degree} of the quotient field of $T$ over the
quotient field of $D$. If $P$ is a prime ideal of $D$, then
$\mathbb{K}(P)$ is denoted the residue field of $D$ in $P$, i.e.,
$D_P/PD_P$, which is canonically isomorphic to the field of
quotients of the integral domain $D/P$.

\begin{defn} Let $D\subseteq T$ be an extension of domain and
$\star$ and $\star'$ are semistar operation on $D$ and $T$
respectively. We say that $D\subseteq T$ satisfies the
\emph{$(\star,\star')$-dimension formula (resp.
$(\star,\star')$-dimension inequality formula)} if for all
$Q\in\QSpec^{\star'}(T)$ such that $(Q\cap D)^{\star}\subsetneq
D^{\star}$, $\star'\text{-}\h(Q)+\td_{\mathbb{K}(Q\cap
D)}(\mathbb{K}(Q))=\star\text{-}\h(Q\cap D)+\td_D(T).$ $
(\text{resp. }\star'\text{-}\h(Q)+\td_{\mathbb{K}(Q\cap
D)}(\mathbb{K}(Q))\leq\star\text{-}\h(Q\cap D)+\td_D(T).)$ The
domain $D$ is said to satisfy the \emph{$\star$-dimension formula
(resp. $\star$-dimension inequality formula)} if for all finitely
generated domain $T$ over $D$, $D\subseteq T$ satisfies the
$(\star,d_T)$-dimension formula (resp. $(\star,d_T)$-dimension
inequality formula).
\end{defn}

If $\star=d_D$ and $\star'=d_T$, then these definitions coincides
with the classical ones (see \cite{Kab}, \cite{AC}, and \cite{FHP}).

\begin{prop}\label{localdiminequality}
Let $D$ be a domain and $\star$ a semistar operation on $D$. Then
the following conditions are equivalent:
\begin{itemize}
\item[(1)] $D$ satisfy
the $\widetilde{\star}$-dimension formula (resp.
$\widetilde{\star}$-dimension inequality formula);

\item[(2)] $D_P$
satisfy the dimension formula for each
$P\in\QSpec^{\widetilde{\star}}(D)$ (resp. dimension inequality
formula);

\item[(3)] $D_M$
satisfy the dimension formula for each
$M\in\QMax^{\widetilde{\star}}(D)$ (resp. dimension inequality
formula).
\end{itemize}
\end{prop}

\begin{proof} We only prove the case of dimension formula and the
other case is the same.

$(1)\Rightarrow(2)$ Let $P\in\QSpec^{\widetilde{\star}}(D)$. Let
$T$ be a finitely generated domain over $D_P$. So that there exist
finitely many elements $\theta_1,\cdots,\theta_r\in T$ such that
$T=D_P[\theta_1,\cdots,\theta_r]$. Set
$T'=D[\theta_1,\cdots,\theta_r]$. Then $T=T'_{D\backslash P}$ and
$T'$ is a finitely generated domain over $D$. Let $Q$ be a prime
ideal of $T$ and set $qD_P:=Q\cap D_P$, where $q(\subseteq P)$ be
a prime ideal of $D$. Thus there exists a prime ideal $Q'$ of
$T'$ such that $Q'\cap(D\backslash P)=\emptyset$ and
$Q=Q'T'_{D\backslash P}$. Thus $Q'\cap D=q$. Since $q\subseteq
P$, we have $q$ is a quasi-$\widetilde{\star}$-prime ideal of
$D$. Since $\widetilde{\star}\text{-}\h(q)=\h(q)$, then by the
hypothesis we have:
$$
\h(Q')+\td_{\mathbb{K}(q)}(\mathbb{K}(Q'))=\h(q)+\td_D(T').
$$
Since $\h(Q')=\h(Q)$ we see that
$$
\h(Q)+\td_{\mathbb{K}(q)}(\mathbb{K}(Q))=\h(q)+\td_{D_P}(T).
$$

$(2)\Rightarrow(3)$ is trivial.

$(3)\Leftarrow(1)$ Suppose that $T$ is a finitely generated domain
over $D$. Let $Q\in\Spec(T)$ and set $P:=Q\cap D$ such that
$P^{\widetilde{\star}}\subsetneq D^{\widetilde{\star}}$. Thus
$P\in\QSpec^{\widetilde{\star}}(D)\cup\{0\}$. Let $M$ be a
quasi-$\widetilde{\star}$-maximal ideal of $D$ containing $P$. Note
that $T_{D\backslash M}$ is a finitely generated domain over $D_M$
and that $Q\cap(D\backslash M)\neq\emptyset$. Thus $QT_{D\backslash
M}\in\Spec(T_{D\backslash M})$ and $PD_M=QT_{D\backslash M}\cap
D_M$. Therefore by the (3), we have
$$
\h(QT_{D\backslash
M})+\td_{\mathbb{K}(PD_M)}(\mathbb{K}(QT_{D\backslash
M}))=\h(PD_M)+\td_{D_M}(T_{D\backslash M}).
$$
Now since $\widetilde{\star}\text{-}\h(Q)=\h(Q)=\h(QT_{D\backslash
M})$, $\h(P)=\h(PD_M)$,
$\td_{\mathbb{K}(P)}(\mathbb{K}(Q))=\td_{\mathbb{K}(PD_M)}(\mathbb{K}(QT_{D\backslash
M}))$ and $\td_D(T)=\td_{D_M}(T_{D\backslash M})$ the proof is
complete.
\end{proof}

Let $D$ be an integral domain with quotient field $K$, let $X$,
$Y$ be two indeterminates over $D$ and let $\star$ be a semistar
operation on D. Set $D_1:=D[X]$, $K_1:=K(X)$ and take the
following subset of $\Spec(D_1)$:
$$\Theta_1^{\star}:=\{Q_1\in\Spec(D_1)|\text{ }Q_1\cap D=(0)\text{ or }(Q_1\cap D)^{\star_f}\subsetneq D^{\star}\}.$$
Set
$\mathfrak{S}_1^{\star}:=\mathcal{S}(\Theta_1^{\star}):=D_1[Y]\backslash(\bigcup\{Q_1[Y]
|Q_1\in\Theta_1^{\star}\})$ and:
$$E^{\circlearrowleft_{\mathfrak{S}_1^{\star}}}:=E[Y]_{\mathfrak{S}_1^{\star}}\cap
K_1, \text{   for all }E\in \overline{\mathcal{F}}(D_1).$$

It is proved in \cite[Theorem 2.1]{S} that the mapping
$\star[X]:=\circlearrowleft_{\mathfrak{S}_1^{\star}}:
\overline{\mathcal{F}}(D_1)\to\overline{\mathcal{F}}(D_1)$,
$E\mapsto E^{\star[X]}$ is a stable semistar operation of finite
type on $D[X]$, i.e., $\widetilde{\star[X]}=\star[X]$. It is also
proved that $\widetilde{\star}[X]=\star_f[X]=\star[X]$,
$d_D[X]=d_{D[X]}$ and
$\QSpec^{\star[X]}(D[X])=\Theta_1^{\star}\backslash\{0\}$. If
$X_1,\cdots,X_r$ are indeterminates over $D$, for $r\geq2$, we let
$$\star[X_1,\cdots,X_r]:=(\star[X_1,\cdots,X_{r-1}])[X_r],$$ where
$\star[X_1,\cdots,X_{r-1}]$ is a stable semistar operation of
finite type on $D[X_1,\cdots,X_{r-1}]$. For an integer $r$, put
$\star[r]$ to denote $\star[X_1,\cdots,X_r]$ and $D[r]$ to denote
$D[X_1,\cdots,X_r]$.

Following \cite{Sah}, the domain $D$ is called
\emph{$\star$-catenary}, if for each pair $P\subset Q$ of
quasi-$\star$-prime ideals  of $D$, any two saturated chain of
quasi-$\star$-prime ideals between $P$ and $Q$ have the same
finite length. If for each $n\geq1$, the polynomial ring $D[n]$ is
$\star[n]$-catenary, then $D$ is said to be
\emph{$\star$-universally catenarian}. Every P$\star$MD which is
$\widetilde{\star}$-LFD (that is $\h(P)<\infty$ for all
$P\in\QSpec^{\widetilde{\star}}(D)$), is
$\widetilde{\star}$-universally catenarian by \cite[Theorem
3.4]{Sah}.

\begin{cor}\label{UCTODIMIN} Let $D$ be an $\widetilde{\star}$-universally catenarian domain. Then $D$ satisfies
the $\widetilde{\star}$-dimension formula.
\end{cor}

\begin{proof} Let $P\in\QSpec^{\widetilde{\star}}(D)$. Hence $D_P$ is a
universally catenarian domain by \cite[Lemma 3.3]{Sah}. Thus by
\cite[Corollary 4.8]{BDF}, $D_P$ satisfies the dimension formula.
Now Proposition \ref{localdiminequality} completes the proof.
\end{proof}

The domain $D$ is called a \emph{$\star$-strong S-domain}, if
each pair of adjacent quasi-$\star$-prime ideals $P_1\subset P_2$
of $D$, extend to a pair of adjacent quasi-$\star[X]$-prime ideals
$P_1[X]\subset P_2[X]$, of $D[X]$. If for each $n\geq1$, the
polynomial ring $D[n]$ is a $\star[n]$-strong S-domain, then $D$
is said to be an \emph{$\star$-stably strong S-domain}. Every
$\widetilde{\star}$-Noetherian,
$\widetilde{\star}$-quasi-Pr\"{u}fer or
$\widetilde{\star}$-universally catenarian domain is
$\widetilde{\star}$-stably strong S-domain by \cite[Corollaries
2.6 and 3.6]{Sah}.

\begin{cor} Let $D$ be an $\widetilde{\star}$-stably strong S-domain. Then $D$ satisfies
the $\widetilde{\star}$-dimension inequality formula.
\end{cor}

\begin{proof} Use \cite[Proposition 2.5]{Sah} and \cite[Theorem 1.6]{Kab} and the same
argument as proof of Corollary \ref{UCTODIMIN}.
\end{proof}

A valuation overring $V$ of $D$ is called a \emph{$\star$-valuation
overring of $D$} provided $F^{\star}\subseteq FV$, for each $F\in
f(D)$. Following \cite{S}, the \emph{$\star$-valuative dimension} of
$D$ is defined as:
$$
\star\text{-}\dim_v(D):=\sup\{\dim(V)|V\text{ is
}\star\text{-valuation overring of }D\}.
$$
Although Example 4.4 of \cite{S} shows that $\star$-$\dim(D)$ is not
always less that or equal to $\star$-$\dim_v(D)$, but it is observed
in \cite{S} that
$\widetilde{\star}$-$\dim(D)\leq\widetilde{\star}$-$\dim_v(D)$. We
say that $D$ is a \emph{$\star$-Jaffard domain}, if
$\star\text{-}\dim(D)=\star\text{-}\dim_v(D)<\infty$. When $\star=d$
the identity operation then $d$-Jaffard domain coincides with the
classical Jaffard domain (cf. \cite{ABDFK}). It is proved in
\cite{S}, that $D$ is a $\widetilde{\star}$-Jaffard domain if and
only if
$$\star[X_1,\cdots,X_n]\text{-}\dim(D[X_1,\cdots,X_n])=\widetilde{\star}\text{-}\dim(D)+n,$$
for each positive integer $n$.

\begin{lem}\label{valuloc} For each domain $D$,
$\widetilde{\star}\text{-}\dim_v(D)=\sup\{\dim_v(D_P)|P\in\QSpec^{\widetilde{\star}}(D)\}$.
\end{lem}

\begin{proof} We can assume that $\widetilde{\star}\text{-}\dim_v(D)$ is a finite number.
Suppose that $n=\widetilde{\star}\text{-}\dim_v(D)$. Then there
exists a $\widetilde{\star}$-valuation overring $V$, with maximal
ideal $N$, of $D$ such that $\dim(V)=n$. Set $P:=N\cap D$. So that
$V$ is a valuation overring of $D_P$. Hence
$n=\dim(V)\leq\dim_v(D_P)\leq\widetilde{\star}\text{-}\dim_v(D)=n$,
where the second inequality is true since each valuation overring of
$D_P$ is a $\widetilde{\star}$-valuation overring of $D$
(\cite[Theorem 3.9]{FL}).
\end{proof}

In \cite[Page 174]{ABDFK} it is proved that a finite-dimensional
domain satisfying the dimension inequality formula is a Jaffard
domain. In the following result we give the semistar analogue of
the mentioned result.

\begin{thm} Let $D$ be a domain of finite $\widetilde{\star}$-dimension. If
$D$ satisfies the $\widetilde{\star}$-dimension inequality
formula, then $D$ is a $\widetilde{\star}$-Jaffard domain.
\end{thm}

\begin{proof} Let $P\in\QSpec^{\widetilde{\star}}(D)$. Then $D_P$ is a finite dimensional domain and
satisfies the dimension inequality formula by Proposition
\ref{localdiminequality}. Consequently $D_P$ is a Jaffard domain
by \cite{ABDFK}. Thus using Lemma \ref{valuloc}, we have
\begin{align*}
\widetilde{\star}\text{-}\dim(D)= & \sup\{\dim(D_P)|P\in\QSpec^{\widetilde{\star}}(D)\} \\[1ex]
                      = & \sup\{\dim_v(D_P)|P\in\QSpec^{\widetilde{\star}}(D)\}\\[1ex]
                      = & \widetilde{\star}\text{-}\dim_v(D).
\end{align*}
Thus $D$ is a $\widetilde{\star}$-Jaffard domain.
\end{proof}

Therefore we have the following implications for finite
$\widetilde{\star}$-dimensional domains:
$$
\begin{array}{ccc}
  \textsf{$\widetilde{\star}$-Noetherian}\text{ or }\textsf{$\widetilde{\star}$-quasi-Pr\"{u}fer} &  & \textsf{P$\star$MD} \\
  \Downarrow &  & \Downarrow \\
  \textsf{$\widetilde{\star}$-stably strong
S-domain} & \Leftarrow &  \textsf{$\widetilde{\star}$-universally
catenary}\\
  \Downarrow &  & \Downarrow \\
  \textsf{$\widetilde{\star}$-dimension inequality formula} & \Leftarrow & \textsf{$\widetilde{\star}$-dimension formula} \\
  \Downarrow &  &  \\
  \textsf{$\widetilde{\star}$-Jaffard} &  &
\end{array}
$$

Let $D$ be a domain with quotient field $K$, let $X$ be an
indeterminate over $D$, let $\star$ be a semistar operation on D,
and let $P$ be a quasi-$\star$-prime ideal of $D$ (or $P=0$). Set
$$
\mathcal{S}_P^{\star}:=(D/P)[X]\backslash\bigcup\{(Q/P)[X] \mid
Q\in\QSpec^{\star_f}(D)\text{ and } P\subseteq Q\}.
$$
Clearly, $\mathcal{S}_P^{\star}$ is a multiplicatively closed
subset of $(D/P)[X]$.

For all $E\in \overline{\mathcal{F}}(D/P)$, set
$$
E^{\circlearrowleft_{\mathcal{S}_P^{\star}}}:=E(D/P)[X]_{\mathcal{S}_P^{\star}}\cap
(D_P/PD_P).
$$
It is proved in \cite[Theorem 3.2]{DS1} that the mapping
$\star/P:=\circlearrowleft_{\mathcal{S}_P^{\star}}:
\overline{\mathcal{F}}(D/P)\to\overline{\mathcal{F}}(D/P)$,
$E\mapsto E^{\circlearrowleft_{\mathcal{S}_P^{\star}}}$, is a
stable semistar operation of finite type on $D/P$; i.e.,
$\widetilde{\star/P}=\star/P$,
$\QMax^{\star/P}(D/P)=\{Q/P\in\Spec(D/P)\mid
Q\in\QMax^{\star_f}(D)\text{ and }P\subseteq Q\}$,
$\widetilde{\star}/P=\star_f/P=\star/P$ and $d_D/P=d_{D/P}$.

\begin{lem}\label{factor} A domain $D$ is $\widetilde{\star}$-universally
catenarian if and only if $D/P$ is $(\star/P)$-universally
catenarian for each $P\in\QSpec^{\widetilde{\star}}(D)$.
\end{lem}

\begin{proof} $(\Rightarrow)$ Let $P\in\QSpec^{\widetilde{\star}}(D)$. By \cite[Theorem 3.2 (a)]{DS1},
$\star/P=\widetilde{\star/P}$. Hence, by \cite[Proposition 3.2
and Lemma 3.3]{Sah}, $D/P$ is $(\star/P)$-universally catenarian
if and only if $(D/P)_{\mathcal{M}}$ is a universally catenarian
domain for each $\mathcal{M} \in \QMax^{\star/P}(D/P)$, that is
(by \cite[Theorem 3.2 (b)]{DS1}), if and only if $D_M/PD_M$ is a
universally catenarian domain whenever $P$ is a subset of $M
\in\QMax^{\widetilde{\star}}(D)$. But by \cite[Proposition 3.2
and Lemma 3.3]{Sah}, $D_M$ is a universally catenarian domain for
all $M \in \QMax^{\widetilde{\star}}(D)$. This, in turn, is
immediate since any factor domain of a universally catenarian
domain must be a universally catenarian domain.

$(\Leftarrow)$ It is enough to consider $P=0$, since we have
$\star/0=\widetilde{\star}$.
\end{proof}

In \cite[Corollary 14.D]{M} it is proved that a Noetherian domain
$D$ is an universally catenarian domain if and only if $D$ is
catenary and $D/P$ satisfies the dimension formula for each
$P\in\Spec(D)$. In the following result we give the semistar
analogue of this result.

\begin{thm} Let $D$ be a $\widetilde{\star}$-Noetherian domain. Then $D$ is an
$\widetilde{\star}$-universally catenarian domain if and only if
$D$ is $\widetilde{\star}$-catenary and $D/P$ satisfies the
$(\star/P)$-dimension formula for each
$P\in\QSpec^{\widetilde{\star}}(D)$.
\end{thm}

\begin{proof} $(\Rightarrow)$ Let $P\in\QSpec^{\widetilde{\star}}(D)$.
Then $D/P$ is $(\star/P)$-universally catenarian by Lemma
\ref{factor}. Hence $D/P$ satisfies the $(\star/P)$-dimension
formula by Corollary \ref{UCTODIMIN}.

$(\Leftarrow)$ Let $M\in\QMax^{\widetilde{\star}}(D)$. It is
enough to show that $D_M$ is a universally catenarian domain. To
this end let $PD_M$ be a prime ideal of $D_M$. Thus $P$ is a
quasi-$\widetilde{\star}$-prime ideal of $D$. Since $D/P$
satisfies the $(\star/P)$-dimension formula, then
$(D/P)_{M/P}=D_M/PD_M$ satisfies the dimension formula by Theorem
\ref{localdiminequality}. On the other hand $D_M$ is a Noetherian
domain by \cite[Proposition 3.8]{EFP} and catenary by
\cite[Proposition 3.2]{Sah}. Consequently $D_M$ is a universally
catenarian domain by \cite[Corollary 14.D]{M}.
\end{proof}

Recall that the celebrated theorem of Ratliff \cite[Theorem
2.6]{Rat} says that a Noetherian ring $R$ is universally
catenarian if and only if $R[X]$ is catenarian. On the other hand
it is proved in \cite[Theorem 1]{BDF1} that the Noetherian
assumption in Ratliff's theorem can be replaced with going-down
condition by proving that: for a going-down domain $D$, we have
$D$ is universally catenarian if and only if $D[X]$ is catenarian
if and only if $D$ is an LFD strong S-domain. As a semistar
analogue in \cite[Theorem 3.7]{Sah} we proved that: suppose that
$D$ is $\widetilde{\star}$-Noetherian. Then $D[X]$ is
$\star[X]$-catenary if and only if $D$ is
$\widetilde{\star}$-universally catenarian. In the last theorem
of this section we treat the second case.

Let $D\subseteq T$ be an extension of domains. Let $\star$ and
$\star'$ be semistar operations on $D$ and $T$, respectively.
Following \cite{DS}, we say that $D\subseteq T$ satisfies
$(\star,\star')$-$\GD$ if, whenever $P_0\subset P$ are
quasi-$\star$-prime ideals of $D$ and $Q$ is a
quasi-$\star'$-prime ideal of $T$ such that $Q\cap D=P$, there
exists a quasi-$\star'$-prime ideal $Q_0$ of $T$ such that
$Q_0\subseteq Q$ and $Q_0\cap D= P_0$. The integral domain $D$ is
said to be a \emph{$\star$-going-down domain} (for short, a
$\star$-$\GD$ {\it domain}) if, for every overring $T$ of $D$ the
extension $D\subseteq T$ satisfies $(\star,d_T)$-$\GD$. These
concepts are the semistar versions of the ``classical'' concepts
of going-down property and the going-down domains (cf.
\cite{DP}). It is known by \cite[Propositions 3.5 and 3.2(e)]{DS}
that every P$\star$MD and every integral domain $D$ with
$\star$-$\dim(D)=1$ is a $\star$-$\GD$ domain.

\begin{thm} Let $D$ be a $\widetilde{\star}$-$\GD$ domain. The
following statements are equivalent:
\begin{itemize}
\item[(1)] $D$ is a $\widetilde{\star}$-LFD $\widetilde{\star}$-strong S-domain.

\item[(2)] $D$ is $\widetilde{\star}$-universally catenarian.

\item[(3)] $D[X]$ is $\star[X]$-catenarian.
\end{itemize}
\end{thm}

\begin{proof} $(1)\Rightarrow(2)$ holds by \cite[Theorem 4.1]{Sah} and $(2)\Rightarrow(3)$ is
trivial. For $(3)\Rightarrow(1)$ let
$P\in\QSpec^{\widetilde{\star}}(D)$. Then $D_P$ is a going-down
domain by \cite[Proposition 2.5]{DS1} and $D_P[X]$ is catenarian
by \cite[Lemma 3.3]{Sah}. Thus $D_P$ is a LFD strong S-domain by
\cite[Theorem 1]{BDF1}. Hence $D$ is a $\widetilde{\star}$-LFD
$\widetilde{\star}$-strong S-domain by \cite[Proposition
2.4]{Sah}.
\end{proof}

\section{Characterizations of $\star$-quasi-Pr\"{u}fer domains}

In this section we give some characterization of
$\widetilde{\star}$-quasi-Pr\"{u}fer domains. We need to recall
the definition of $(\star,\star')$-linked overrings. Let $D$ be a
domain and $T$ an overring of $D$. Let $\star$ and $\star'$ be
semistar operations on $D$ and $T$, respectively. One says that
$T$ is \emph{$(\star,\star')$-linked to} $D$ (or that $T$ is a
$(\star,\star')${\it -linked overring of} $D$) if
$F^{\star}=D^{\star}\Rightarrow (FT)^{\star'}=T^{\star'}$, when
$F$ is a nonzero finitely generated ideal of $D$ (cf. \cite{EF}).
In particular we are interested in the case $\star'=d_T$. We
first recall the following characterization of
$\widetilde{\star}$-quasi-Pr\"{u}fer domains.

\begin{thm}(\cite[Theorem 4.3]{Sah})\label{qJaf} Let $D$ be an integral domain.
Suppose that $\widetilde{\star}$-$\dim(D)$ is finite. Consider
the following statements:
\begin{itemize}
\item[(1')] Each $(\star,\star')$-linked overring $T$ of $D$ is an
$\widetilde{\star'}$-universally catenarian domain.

\item[(1)] Each $(\star,\star')$-linked overring $T$ of $D$ is an
$\widetilde{\star'}$-stably strong S-domain.

\item[(2)] Each $(\star,\star')$-linked overring $T$ of $D$ is an
$\widetilde{\star'}$-strong S-domain.

\item[(3)] Each $(\star,\star')$-linked overring $T$ of $D$ is an
$\widetilde{\star'}$-Jaffard domain.

\item[(4)] Each $(\star,\star')$-linked overring $T$ of $D$ is an
$\widetilde{\star'}$-quasi-Pr\"{u}fer domain.

\item[(5)] $D$ is an $\widetilde{\star}$-quasi-Pr\"{u}fer domain.
\end{itemize}
Then
$(1')\Rightarrow(1)\Leftrightarrow(2)\Leftrightarrow(3)\Leftrightarrow(4)\Leftrightarrow(5)$.
\end{thm}

\begin{proof} The implication $(1')\Rightarrow(1)$ holds by \cite[Corollary
3.6]{Sah} and $(1)\Rightarrow(2)$ is trivial. For
$(2)\Rightarrow(5)$ see proof of \cite[Theorem 4.3]{Sah} part
$(3)\Rightarrow(6)$. The implication $(5)\Rightarrow(1)$ holds by
\cite[Corollary 2.6]{Sah}. For
$(4)\Leftrightarrow(5)\Leftrightarrow(6)$ see \cite[Theorem
4.14]{S}.
\end{proof}

Now we have the following theorem; a result reminiscent of the
well-known result of Ayache, Cahen and Echi \cite{ACE} (see also
\cite[Theorem 6.7.8]{FHP}) for quasi-Pr\"{u}fer domains.

\begin{thm} Let $D$ be an integral domain. Suppose that $\widetilde{\star}$-$\dim(D)$ is finite. Then the
following statements are equivalent:
\begin{itemize}
\item[(1)] Each $(\star,d_T)$-linked overring $T$ of $D$ is a
stably strong S-domain.

\item[(2)] Each $(\star,d_T)$-linked overring $T$ of $D$ is a
strong S-domain.

\item[(3)] Each $(\star,d_T)$-linked overring $T$ of $D$ is a
Jaffard domain.

\item[(4)] Each $(\star,d_T)$-linked overring $T$ of $D$ is a
quasi-Pr\"{u}fer domain.

\item[(5)] $D$ is an $\widetilde{\star}$-quasi-Pr\"{u}fer domain.
\end{itemize}
\end{thm}

\begin{proof} We only prove the equivalence of $(1)\Leftrightarrow(5)$ and the proofs of
$(2)\Leftrightarrow(5)$, $(3)\Leftrightarrow(5)$, and
$(4)\Leftrightarrow(5)$, are the same. The implication
$(5)\Rightarrow(1)$ holds by Theorem \ref{qJaf}. For
$(1)\Rightarrow(5)$ let $P\in\QSpec^{\widetilde{\star}}(D)$. It
is enough for us to show that $D_P$ is a quasi-Pr\"{u}fer domain
by \cite[Theorem 2.16]{CF}. To this end let $T$ be an overring of
$D_P$. Then $T_{D\backslash P}=T$ and therefore $T$ is
$(\star,d_T)$-linked overring of $D$ by \cite[Example 3.4
(1)]{EF}. Thus by the hypothesis we have $T$ is a stably strong
S-domain. Therefore $D_P$ is a quasi-Pr\"{u}fer domain by
\cite[Theorem 6.7.8]{FHP}.
\end{proof}

\begin{thm}\label{qdif} Let $D$ be an integral domain. Suppose that $\widetilde{\star}$-$\dim(D)$ is finite.
Then the following statements are equivalent:
\begin{itemize}
\item[(1)] $D$ is an $\widetilde{\star}$-quasi-Pr\"{u}fer domain.

\item[(2)] For each $(\star,\star')$-linked overring $T$ of $D$, every extension of domains $T\subseteq S$, satisfies the
$(\widetilde{\star'},\widetilde{\star''})$-dimension inequality
formula, where $\star'$ and $\star''$ are semistar operations on
$T$ and $S$ respectively.
\end{itemize}
\end{thm}

\begin{proof} $(1)\Rightarrow(2)$ If $D$ is an $\widetilde{\star}$-quasi-Pr\"{u}fer
domain and $T$ is $(\star,\star')$-linked to $D$, then $T$ is a
$\widetilde{\star'}$-Jaffard domain by \cite[Theorem 4.14]{S}. Let
$Q\in\QSpec^{\widetilde{\star''}}(S)$ such that $(Q\cap
T)^{\widetilde{\star'}}\subsetneq T^{\widetilde{\star'}}$ and set
$q:=Q\cap T$. Then we have
$q\in\QSpec^{\widetilde{\star'}}(T)\cup\{0\}$. Set $P:=q\cap D$.
Thus we have $P\in\QSpec^{\widetilde{\star}}(D)\cup\{0\}$. Therefore
$D_P$ and hence $T_q$, are quasi-Pr\"{u}fer domains by \cite[Theorem
1.1]{CF}. In particular $T_q$ is a Jaffard domain. So that we have
\begin{align*}
\dim(S_Q)+\td_{\mathbb{K}(q)}(\mathbb{K}(Q))\leq & \dim_v(S_Q)+\td_{\mathbb{K}(q)}(\mathbb{K}(Q)) \\[1ex]
                      \leq & \dim_v(T_q)+\td_{T}(S),
\end{align*}
where the first inequality holds since $\dim(S_Q)\leq\dim_v(S_Q)$
and the second one is by \cite[Lemma 6.7.3]{FHP}. The conclusion
follows easily from the fact that $\dim(T_q)=\dim_v(T_q)$.

$(2)\Rightarrow(1)$ Let $T$ be an overring of $D$ and $\star'$ be
a semistar operation on $T$ such that $T$ is
$(\star,\star')$-linked to $D$. Let $(V,N)$ be any
$\widetilde{\star'}$-valuation overring of $T$. Then $V$ is
$(\star',d_V)$-linked to $T$ by \cite[Lemma 2.7]{EFP}. Set
$Q:=N\cap T$. Then by assumption we have
$$\dim(V)\leq\dim(T_Q)-\td_{\mathbb{K}(Q)}(\mathbb{K}(N)).$$
In particular
$\dim(V)\leq\dim(T_Q)\leq\widetilde{\star'}\text{-}\dim(T)$, and
hence
$\widetilde{\star'}\text{-}\dim_v(T)=\widetilde{\star'}\text{-}\dim(T)$,
that is $T$ is a $\widetilde{\star'}$-Jaffard domain. Thus $D$ is an
$\widetilde{\star}$-quasi-Pr\"{u}fer domain by \cite[Theorem
4.14]{S}.
\end{proof}

\begin{cor}\label{d} Let $D$ be an integral domain. Suppose that $\widetilde{\star}$-$\dim(D)$ is finite. Then the
following statements are equivalent:
\begin{itemize}
\item[(1)] $D$ is an $\widetilde{\star}$-quasi-Pr\"{u}fer domain.

\item[(2)] For each $(\star,d_T)$-linked overring $T$ of $D$, every extension of domains $T\subseteq S$, satisfies the
dimension inequality formula.
\end{itemize}
\end{cor}

\begin{proof} $(1)\Rightarrow(2)$ holds by Theorem \ref{qdif}. For
$(2)\Rightarrow(1)$ let $P\in\QSpec^{\widetilde{\star}}(D)$. It
is enough for us to show that $D_P$ is a quasi-Pr\"{u}fer domain
by \cite[Theorem 2.16]{CF}. To this end let $T$ be an overring of
$D_P$. Then $T_{D\backslash P}=T$ and therefore $T$ is
$(\star,d_T)$-linked overring of $D$ by \cite[Example 3.4
(1)]{EF}. If $T\subseteq S$ is any extension of domains, then
$T\subseteq S$ satisfies the dimension inequality formula by the
hypothesis. Therefore $D_P$ is a quasi-Pr\"{u}fer domain by
\cite[Theorem 6.7.4]{FHP}.
\end{proof}

Recall that an integral domain $D$ is called a \emph{UM$t$-domain}
if every upper to zero in $D[X]$ is a maximal $t$-ideal
 and has been studied by several authors (See
\cite{HZ}, \cite{FGH} and \cite{CF}). It is observed in
\cite[Corollary 2.4 (b)]{CF} that $D$ is a $w$-quasi-Pr\"{u}fer
domain if and only if $D$ is a UM$t$-domain. The following corollary
is a new characterization of UM$t$ domains.

\begin{cor} Let $D$ be an integral domain. Suppose that $w$-$\dim(D)$ is finite. Then the
following statements are equivalent:
\begin{itemize}
\item[(1)] Each $(t_D,d_T)$-linked overring $T$ of $D$ is a
stably strong S-domain.

\item[(2)] Each $(t_D,d_T)$-linked overring $T$ of $D$ is a
strong S-domain.

\item[(3)] Each $(t_D,d_T)$-linked overring $T$ of $D$ is a
Jaffard domain.

\item[(4)] Each $(t_D,d_T)$-linked overring $T$ of $D$ is a
quasi-Pr\"{u}fer domain.

\item[(5)] For each $(t_D,d_T)$-linked overring $T$ of $D$, every extension of domains $T\subseteq S$, satisfies the
dimension inequality formula.

\item[(6)] $D$ is a UM$t$ domain.
\end{itemize}
\end{cor}

\section{Arnold's formula}

In the last section we extends some results of J. Arnold of the
dimension of polynomial rings to the setting of the semistar
operations. First we wish to give the following lemma which is a new
property of semistar valuative dimension.

\begin{lem}\label{Jaffard}(see \cite[Theorem 4.2]{S}) Let $D$ be an integral domain
and $n$ be an integer. Then the following statements are
equivalent:
\begin{itemize}
\item[(1)] Each $(\star,d_T)$-linked overring $T$ of $D$ has
dimension at most $n$.

\item[(2)] Each $\widetilde{\star}$-valuation overring of $D$ has
dimension at most $n$.
\end{itemize}
\end{lem}

\begin{proof} The implication $(1)\Rightarrow(2)$ is trivial. For
$(2)\Rightarrow(1)$ let $T$ be a $(\star,d_T)$-linked overring of
$D$ and $V$ be a valuation overring of $T$. Then it is easy to
see that $V$ is $(\star,d_V)$-linked overring of $D$. Thus by
\cite[Lemma 2.7]{EFP}, $V$ is an $\widetilde{\star}$-valuation
overring of $D$. Hence $\dim(V)\leq n$. Consequently
$\dim(T)\leq\dim_v(T)\leq n$ as desired.
\end{proof}

When $\star=d_D$, the equivalence of (1) and (3) of the following
theorem is due to J. Arnold \cite[Theorem 6]{A}.

\begin{thm} Let $D$ be an integral domain, and $n$ be an integer. Then the following statements are
equivalent:
\begin{itemize}
\item[(1)] $\widetilde{\star}$-$\dim_v(D)=n$.

\item[(2)] $\star[n]\text{-}\dim(D[n])=2n$.

\item[(3)] $\star[r]\text{-}\dim(D[r])=r+n$ for all $r\geq n-1$.

\item[(4)] Each $(\star,d_T)$-linked overring $T$ of $D$ has
dimension at most $n$, and $n$ is minimal.
\end{itemize}
\end{thm}

\begin{proof} The equivalence $(1)\Leftrightarrow(2)$ follows from
\cite[Theorem 4.5]{S}, and $(3)\Rightarrow(2)$ is trivial. For
$(1)\Rightarrow(3)$ suppose that
$\widetilde{\star}$-$\dim_v(D)=n$. Then For all $r\geq n$ we have
$\star[r]\text{-}\dim(D[r])=\star[r]\text{-}\dim_v(D[r])=r+\widetilde{\star}\text{-}\dim_v(D)=r+n$,
by \cite[Corollary 4.7 and Theorem 4.8]{S}. Now assume that
$r=n-1$. Since $\widetilde{\star}$-$\dim_v(D)=n$, there exists a
quasi-$\widetilde{\star}$-prime ideal $M$ of $D$ such that
$n=\dim_v(D_M)$, by Lemma \ref{valuloc}. So that by \cite[Theorem
6]{A} we have $\dim(D_M[r])=r+n$. Let
$\mathcal{P}\in\QSpec^{\star[r]}(D[r])$ be such that
$\star[r]\text{-}\dim(D[r])=\h(\mathcal{P})$. Set
$P:=\mathcal{P}\cap D$. Then by \cite[Remark 2.3]{S} we have
$P\in\QSpec^{\widetilde{\star}}(D)\cup\{(0)\}$. Thus
\begin{align*}
r+n\leq & \star[r]\text{-}\dim(D[r])=\h(\mathcal{P}) \\[1ex]
   = & \dim(D[r]_{\mathcal{P}})=\dim(D_P[r]_{\mathcal{P}D_P[r]})\\[1ex]
   \leq & \dim(D_P[r])\leq\dim(D_M[r])=r+n,
\end{align*}
where the first inequality holds by \cite[Theorem 3.1]{S}. Hence
$\star[r]\text{-}\dim(D[r])=r+n$ for all $r\geq n-1$.

The equivalence $(1)\Leftrightarrow(4)$ follows from Lemma
\ref{Jaffard}.
\end{proof}

As an immediate consequence we have:

\begin{cor}\label{loc} $\widetilde{\star}\text{-}\dim_v(D)=\sup\{\dim(T)| T\text{ is
}(\star,d_T)\text{-linked overring of }D\}$.
\end{cor}

One of the famous formulas in the dimension theory of commutative
rings is the Arnold's formula \cite[Theorem 5]{A} which states as
$$
\dim(D[n])=n+\sup\{\dim(D[\theta_1,\cdots,\theta_n])|\{\theta_i\}_1^n\subseteq
K\}.
$$
Now we prove the semistar analogue of Arnold's formula.

\begin{lem}\label{AG} Let $D$ be an integral domain
and $n$ be an integer. Then
$$
\star[n]\text{-}\dim(D[n])=\sup\{\dim(D_M[n])|M\in\QMax^{\widetilde{\star}}(D)\}.
$$
\end{lem}

\begin{proof} If $P$ is a quasi-$\widetilde{\star}$-prime ideal of $D$, and if $QD_P[n]$ is
a non-zero prime ideal of $D_P[n](=D[n]_{D\backslash P})$, then
$Q\cap D\subseteq P$ and hence $Q\in\QSpec^{\star[n]}(D[n])$ by
\cite[Remark 2.3]{S}. So that the inequality $\geq$ is true. Now
let $Q\in\QMax^{\star[n]}(D[n])$ be such that
$\star[n]\text{-}\dim(D[n])=\h(Q)$, and set $P:=Q\cap D$. Then by
\cite[Remark 2.3]{S} we have
$P\in\QSpec^{\widetilde{\star}}(D)\cup\{0\}$. So that
\begin{align*}
\star[n]\text{-}\dim(D[n])= & \h(Q)=\dim(D[n]_Q) \\[1ex]
   = & \dim(D_P[n]_{QD_P[n]})\leq\dim(D_P[n]) \\[1ex]
   \leq & \star[n]\text{-}\dim(D[n]).
\end{align*}
Therefore the proof is complete.
\end{proof}

\begin{cor} Let $D$ be an integral domain
and $n$ be an integer. Then there exist a
quasi-$\widetilde{\star}$-maximal ideal $M$ of $D$ and a
quasi-$\star[n]$-maximal ideal $Q$ of $D[n]$ such that $M=Q\cap D$
and
$$
\star[n]\text{-}\dim(D[n])=\h(Q)=n+\h(M[n]).
$$
\end{cor}

\begin{proof} By Lemma \ref{AG} there exists a quasi-$\widetilde{\star}$-maximal ideal $M$ of
$D$ such that $\star[n]\text{-}\dim(D[n])=\dim(D_M[n])$. Thus
there exists a prime ideal $Q$ of $D[n]$ such that
$Q\cap(D\backslash M)=\emptyset$, $\dim(D_M[n])=\h(QD_M[n])$ and
that $QD_M[n]$ is a maximal ideal of $D_M[n]$. Since $Q\cap
D\subseteq M$ we have $Q$ is a quasi-$\star[n]$-prime of $D[n]$
by \cite[Remark 2.3]{S}, and since
$\star[n]\text{-}\dim(D[n])=\h(Q)$, we have $Q$ is a
quasi-$\star[n]$-maximal ideal of $D[n]$. Set $PD_M:=QD_M[n]\cap
D_M$ for some $P\in\QSpec^{\widetilde{\star}}(D)$. Then by
\cite[Corollary 2.9]{AG} we have $\h(QD_M[n])=n+\h(PD_M[n])$ and
that $PD_M$ is a maximal ideal of $D_M$. Thus we have $P=M$ and
$$
\star[n]\text{-}\dim(D[n])=\h(QD_M[n])=n+\h(MD_M[n])=n+\h(M[n]),
$$
which ends the proof.
\end{proof}

We are ready to prove the semistar analogue of Arnold's formula.

\begin{thm}\label{A} Let $D$ be an integral domain
and $n$ be a positive integer. Then
$$
\star[n]\text{-}\dim(D[n])=n+\sup\{\widetilde{\star}_{\iota}\text{-}\dim(D[\theta_1,\cdots,\theta_n])|
\{\theta_i\}_1^n\subseteq
K\}.
$$
where $\iota$ is the inclusion map of $D$ in
$D[\theta_1,\cdots,\theta_n]$.
\end{thm}

\begin{proof} Let $M\in\QMax^{\widetilde{\star}}(D)$ and
$\{\theta_i\}_1^n\subseteq K$. Let $Q$ be a maximal ideal of
$D_M[\theta_1,\cdots,\theta_n]$ such that
$\dim(D_M[\theta_1,\cdots,\theta_n])=\h(Q)$. Let $Q_0$ be a prime
ideal of $D[\theta_1,\cdots,\theta_n]$ such that
$Q_0\cap(D\backslash M)=\emptyset$ and
$Q=Q_0D_M[\theta_1,\cdots,\theta_n]$. Thus $Q_0$ is a
quasi-$\widetilde{\star}_{\iota}$-prime ideal of
$D[\theta_1,\cdots,\theta_n]$ since $Q_0\cap D\subseteq M$
(\cite[Remark 2.3]{S}). Hence we obtain that
$\dim(D_M[\theta_1,\cdots,\theta_n])=\h(Q)=\h(Q_0)\leq
\widetilde{\star}_{\iota}\text{-}\dim(D[\theta_1,\cdots,\theta_n])$.
Using Lemma \ref{AG} and Arnold's formula \cite[Theorem 5]{A}, we
have
$$\star[n]\text{-}\dim(D[n])=n+\sup\{\dim(D_M[\theta_1,\cdots,\theta_n])\},$$
where the supremum is taken over
$M\in\QMax^{\widetilde{\star}}(D)$ and $\{\theta_i\}_1^n\subseteq
K$. So that $\star[n]\text{-}\dim(D[n])\leq
n+\sup\{\widetilde{\star}_{\iota}\text{-}\dim(D[\theta_1,\cdots,\theta_n])|
\{\theta_i\}_1^n\subseteq K\}$. Now choose
$M\in\QMax^{\widetilde{\star}}(D)$ and $\{\theta_i\}_1^n\subseteq
K$ such that
$\star[n]\text{-}\dim(D[n])=n+\dim(D_M[\theta_1,\cdots,\theta_n])$.
Let $Q'$ be a quasi-$\widetilde{\star}_{\iota}$-prime ideal of
$D[\theta_1,\cdots,\theta_n]$ such that
$\widetilde{\star}_{\iota}\text{-}\dim(D[\theta_1,\cdots,\theta_n])=\h(Q')$
and set $P':=Q'\cap D$. Thus
\begin{align*}
\widetilde{\star}_{\iota}\text{-}\dim(D[\theta_1,\cdots,\theta_n])= &\h(Q')=\dim(D[\theta_1,\cdots,\theta_n]_{Q'}) \\[1ex]
   = & \dim(D_{P'}[\theta_1,\cdots,\theta_n]_{Q'D_{P'}[\theta_1,\cdots,\theta_n]})\\[1ex]
   \leq &
   \dim(D_{P'}[\theta_1,\cdots,\theta_n])\leq\dim(D_M[\theta_1,\cdots,\theta_n]).
\end{align*}
Hence by the first part of the proof
$\dim(D_M[\theta_1,\cdots,\theta_n])=
\widetilde{\star}_{\iota}\text{-}\dim(D[\theta_1,\cdots,\theta_n])$.
Thus we have
$\star[n]\text{-}\dim(D[n])=n+\widetilde{\star}_{\iota}\text{-}\dim(D[\theta_1,\cdots,\theta_n])$
to complete the proof.
\end{proof}



\end{document}